\documentclass{article}
\usepackage{amsmath,amssymb,latexsym,amsthm,amscd}
\font\goth=eusm10

\newcommand\bc{\mathbf C}

\newcommand\E{\mathcal E}

\newcommand\bz{\mathbf{Z}}
\newcommand\br{\mathbf{R}}
\newcommand\bq{\mathbf{Q}}

\newcommand\bn{\mathbf{N}}

\newcommand\Oc{\hbox{\goth O}}

\newcommand\wH{\widetilde{H}}
\newcommand\wth{\widetilde{h}}

\newcommand\rk{\text{rk\ }}

\newcommand\Aut{\text{Aut}\,}

\newtheorem{proposition}{Proposition}
\newtheorem{theorem}{Theorem}

\newtheorem{definition}{Definition}
\newtheorem{corollary}{Corollary}
\newtheorem{remark}{Remark}
\newtheorem{lemma}{Lemma}

\numberwithin{proposition}{section}
\numberwithin{definition}{section}
\numberwithin{corollary}{section} \numberwithin{remark}{section}
\numberwithin{lemma}{section} \numberwithin{equation}{section}

\numberwithin{theorem}{section}
\newtheorem{example}{Example}

\numberwithin{question}{section} \numberwithin{case}{section}
\numberwithin{example}{section} \numberwithin{conjecture}{section}
\newtheorem{problem}{Problem}[section]

\begin{document}
\title{Correspondences of a K3 surface with itself via moduli of sheaves. I}
\author{Viacheslav V. Nikulin \footnote{Supported by EPSRC grant
EP/D061997/1}} \maketitle

\begin{abstract} Let $X$ be an algebraic K3 surface, $v=(r,H,s)$ a primitive
isotropic Mukai vector on $X$ and $M_X(v)$ the moduli of sheaves
over $X$ with $v$. Let $N(X)$ be Picard lattice of $X$.

In \cite{Nik5} and \cite{Mad-Nik5}, all divisors in moduli of
$(X,H)$ (i. e. pairs $H\in N(X)$ with  $\rk N(X)=2$) implying
$M_X(v)\cong X$ were described. They give some Mukai's
correspondences of $X$ with itself.

Applying these results, we show that there exists $v$ and a
codimension $2$ submoduli in moduli of $(X,H)$ (i. e. a pair $H\in
N(X)$ with $\rk N(X)=3$) implying $M_X(v)\cong X$, but this
submoduli cannot be extended to a divisor in moduli with the same
property. There are plenty of similar examples.

We discuss the general problem of description of all  similar
submoduli (to generalize results of \cite{Nik5} and
\cite{Mad-Nik5} to the general case), and defined by them Mukai's
correspondences of $X$ with itself and their compositions, trying
to outline a possible general theory.
\end{abstract}

\section{Introduction} \label{introduction}

We consider algebraic K3 surfaces $X$ over $\bc$. We remind that a
non-singular projective algebraic (or K\"ahlerian compact) surface
$X$ is a K3 surface if its canonical class $K_X$ is equal to zero
and the irregularity $q=\dim \Omega^1[X]=0$.

Further, $N(X)$ denotes the Picard lattice of $X$, and $T(X)$ the
transcendental lattice. Further, $\rho (X)=\rk N(X)$ denotes the
Picard number of $X$.

We consider primitive isotropic Mukai vectors
\begin{equation}
v=(r,l,s),\ r\in \bn,\ s\in \bz,\  l\in N(X),\  l^2=2rs.
\label{mukaivectorint}
\end{equation}
on $X$. We denote by $Y=M_X(v)=M_X(r,l,s)$ the K3 surface which is the
minimal resolution of singularities of the moduli space of sheaves
over $X$ with the Mukai vector $v$. Clarify in Mukai \cite{Muk1}
--- \cite{Muk4} and Yoshioka \cite{Yoshioka1}. In this case, the corresponding
quasi-universal sheave on $X\times Y$, and its Chern class define
a 2-dimensional algebraic cycle on $X\times Y$ and a
correspondence between $X$ and $Y$ according to Mukai
\cite{Muk2}. It has very nice geometric properties. For more
details, see Sect. \ref{seccompositions}.

If $Y\cong X$, then we obtain an important $2$-dimensional
algebraic cycle on $X\times X$, and a correspondence of $X$ with
itself. Thus, the question, when $Y\cong X$ is very interesting.

If $\rho (X)=1$,  we give the corresponding results in  Sect.
\ref{secpicard1}. The author believes that they should be known to
specialists.

Let $D\in N(X)$.
Then one has the natural isomorphism given by the tensor product
$$
T_D:M_X(r,l,s)\cong M_X(r,l+rD,s+r(D^2/2)+D\cdot l),\ \
\E\mapsto \E\otimes \Oc(D).
$$

For $r,\,s\,>0$, one has the isomorphism which is called {\it reflection}
$$
\delta:M_X(r,l,s)\cong M_X(s,l,r).
$$
E. g. see \cite{Muk2} and \cite{Tyurin1}, \cite{Tyurin2},
\cite{Yoshioka2}.

For integers $d_1,\,d_2>0$ such that
$(d_1,d_2)=(d_1,s)=(r,d_2)=1$, one has the isomorphism
$$
\nu(d_1,d_2):M_X(r,l,s)\cong M_X(rd_1^2,ld_1d_2,sd_2^2)
$$
and its inverse $\nu(d_1,d_2)^{-1}$. See \cite{Muk2}, \cite{Muk3},
\cite{Nik3}, \cite{Mad-Nik5}.

In Theorem \ref{theoremunivred1} and Corollary
\ref{corollaryuniv}, we show that if $\rho (X)=1$ and $X$ is
general, then for two primitive isotropic Mukai vectors $v_1$ and
$v_2$, moduli $M_X(v_1)$ and $M_X(v_2)$ are isomorphic if and only
if there exists their isomorphism which is a composition of the
three isomorphisms described above. They and their specializations
for higher Picard numbers give {\it universal isomorphisms}
between moduli of sheaves over $X$.

It is known
(e.g. see \cite{Tyurin2}) that for $l\in N(X)$ and $\pm l^2>0$, one
has Tyurin's isomorphism
\begin{equation}
\label{Tyurinisomint}
Tyu=Tyu(\pm l)\,:\, M_X(\pm l^2/2,l,\pm 1)\cong X \,.
\end{equation}

In Corollary \ref{corollaryunivred1congX} (see also Remark
\ref{remarkrho=1}) we show that, if $\rho (X)=1$, then
$M_X(r,H,s)$ and $X$ are isomorphic if and only if there exists
their isomorphism which is a  composition of the three universal
isomorphisms above between moduli of sheaves, and the Tyurin's
isomorphism. See also \cite{Nik5} for a similar result.

It was shown in \cite{Nik5} and geometrically interpreted in
\cite{Mad-Nik5} (together with Carlo Madonna) that analogous
results are valid for $\rho (X)=2$ and $X$ which is general for
its Picard lattice, i. e. the automorphism group of the
transcendental periods is trivial: $\Aut (T(X), H^{2,0}(X))=\pm 1$
(see important particular cases of these results in
\cite{Mad-Nik1}, \cite{Mad-Nik2}, \cite{Nik4} and \cite{Mad-Nik3},
\cite{Mad-Nik4}). We review these results in Sect.
\ref{secpicard2}. See Theorems \ref{maintheorem},
\ref{maintheorem2} and \ref{maintheorem3} for exact formulations.
They show that in this case (i. e. when $\rho (X)=2$ and $X$ is
general for its Picard lattice), $M_X(r,H,s)\cong X$  if and only
if there exists an isomorphism between $M_X(r,H,s)$ and $X$ which
is given by a composition of the universal isomorphisms $T_D$,
$\delta$ and $\nu(d_1,d_2)$ between moduli of sheaves over $X$ and
by Tyurin's isomorphism between moduli of sheaves over $X$ and $X$
itself. The results for $\rho (X)=1$ above clarify appearance of
the natural isomorphisms $T_D$, $\delta$, $\nu(d_1,d_2)$, $Tyu$ in
these results for Picard number $2$

Importance of the results for $\rho (X)=2$ and general $X$ is that
they describe all divisorial conditions on moduli of algebraic
polarized K3 surfaces $(X,H)$ which imply that $M_X(r,H,s)\cong
X$. More exactly, the results for $\rho (X)=2$ describe all
abstract polarized Picard lattices $H\in N$ of the rank $\rk N=2$
such that if $H\in N\subset N(X)$, then $M_X(r,H,s)\cong X$. We
remind that the codimension of moduli of such $X$ in the
19-dimensional moduli of algebraic polarized K3 surfaces is then
$1$.  Applying these results, in  Theorems \ref{theoremnecessary}
and \ref{divconditionsnes}, we give a necessary condition on Mukai
vector $(r,H,s)$ and a K3 surface $X$, for the isomorphism
$M_X(r,H,s)\cong X$ would follow from a divisorial condition on
moduli of polarized K3 surfaces. In Example \ref{example5,13},  we
give an exact numerical example when this necessary condition is
not satisfied. Thus, for K3 surfaces $X$ of this example, the
isomorphism $M_X(r,H,s)\cong X$ cannot follow from any divisorial
condition on moduli of polarized K3 surfaces which implies the
isomorphism $M_{X^\prime}(r,H,s)\cong X^\prime$ for K3 surfaces
$X^\prime$ satisfying this divisorial condition.

Applying these results, in Sect. \ref{secpicard3}, Theorem
\ref{theoremrho=3}, we give an exact example of the type of a
primitive isotropic Mukai vector $(r,H,s)$ and a pair $H\in N$ of
a (abstract) polarized K3 Picard lattice or the rank $\rk N=3$
such that for any polarized K3 surface $(X,H)$ with $H\in N\subset
N(X)$ one has $M_X(r,H,s)\cong X$, but this isomorphism does not
follow from any divisorial condition (i. e. from Picard number 2)
on moduli of polarized K3 surfaces. Thus, moduli of these
polarized K3 surfaces have codimension 2, and they cannot be
extended to a divisor in moduli of polarized K3  surfaces
preserving the isomorphism $M_X(r,H,s)\cong X$. This is {\it the
main result of this paper.} In this section, we give many similar
examples for Picard number $\rho (X)\ge 3$.

These results give important corollaries for higher Picard number
$\rho (X)\ge 3$ of the results for Picard number 1 and 2,
described above. They also show that the case of $\rho (X)\ge 3$
is very non-trivial. These are the main subjects of this paper.

Another important subject of this paper is to formulate some
general concepts and predict the general structure of possible
results for higher Picard number $\rho (X)\ge 3$.

At the end of Sect. \ref{secpicard3}, for a type $(r,H,s)$ of a
primitive isotropic Mukai vector,  we introduce the concept of a
{\it critical polarized K3 Picard lattice} $H\in N$ (for the
problem of correspondences of a K3 surface with itself). Roughly
speaking, it means that $M_X(r,H,s)\cong X$ for a polarized K3
surface $X$ such that $H\in N\subset N(X)$, but the same is not
valid for any primitive sublattice $H\in N_1\subset N$ of strictly
smaller rank. Thus, the corresponding moduli of K3 surfaces have
dimension $20-\rk N$, and they are not specialization of moduli of
higher dimension of analogous K3 surfaces.

{\it Classification of critical polarized K3 Picard lattices} is
the main problem of correspondences of a K3 surface with itself
via moduli of sheaves. Our results for $\rho=1$ and $\rho=2$ can
be interpreted as classification of all critical polarized K3
Picard lattices of the rank one and two. The described above
example of the lattice $N$ of the rank three from Theorem
\ref{theoremrho=3}, gives an example of a critical polarized K3
Picard lattice of the rank three. In Theorem \ref{theoremrankcrit}
we prove that $\rk N\le 12$ for a critical polarized K3 Picard
lattice $N$. In Problem \ref{problemcritrank},  we raise a problem
of exact estimate of the rank of critical polarized K3 Picard
lattices for the fixed type of  a primitive isotropic Mukai
vector. This problem is now solved only for very special types,
when the rank is equal to one.

\medskip

In Sect. \ref{seccompositions}, we interpret the results above in
terms of the action of correspondences as isometries in
$H^2(X,\bq)$, and their compositions. For example, Tyurin's
isomorphisms  \ref{Tyurinisomint} give reflections in elements
$l\in N(X)$ and generate the full automorphism group
$O(N(X)\otimes \bq)$. Each isotropic primitive Mukai vector
$(r,H,s)$ on $X$ and isomorphism $M_X(r,H,s)\cong X$ generates
then some class of isometries from $O(N(X)\otimes \bq)$. See Sect.
\ref{seccompositions} for exact formulations. Thus, the main
problem of correspondences of $X$ with itself via moduli of
sheaves is to find all these additional generators and their
relations. We formulate the corresponding problems (1), (2), (3)
and (4) at the end of Sect. \ref{seccompositions}. They show that,
in principle, general results for any $\rho (X)$ should look
similar to the now known results for $\rho (X)=1,\,2$.

\medskip

A reader can see that our general idea is that a very complicated
structure of correspondences of a general (for its Picard lattice)
K3 surface $X$ with itself via moduli of sheaves is hidden inside
of the abstract Picard lattice $N(X)$, and we try to recover this
structure. This should lead to some non-trivial constructions
related to the abstract Picard lattice $N(X)$ and more closely
relate it to geometry of the K3 surface.

\medskip

The author is grateful to D.O. Orlov for useful discussions.

\medskip

\section{Preliminary notations about lattices}
\label{head1}

We use notations and terminology from \cite{Nik2} about lattices,
their discriminant groups and forms. A {\it lattice} $L$ is a
non-degenerate integral symmetric bilinear form. I. e. $L$ is a
free $\bz$-module equipped with a symmetric pairing $x\cdot y\in
\bz$ for $x,\,y\in L$, and this pairing should be non-degenerate.
We denote $x^2=x\cdot x$. The {\it signature} of $L$ is the
signature of the corresponding real form $L\otimes \br$. The
lattice $L$ is called {\it even} if $x^2$ is even for any $x\in
L$. Otherwise, $L$ is called {\it odd}. The {\it determinant} of
$L$ is defined to be $\det L=\det(e_i\cdot e_j)$ where $\{e_i\}$
is some basis of $L$. The lattice $L$ is {\it unimodular} if $\det
L=\pm 1$. The {\it dual lattice} of $L$ is
$L^\ast=Hom(L,\,\bz)\subset L\otimes \bq$. The {\it discriminant
group} of $L$ is $A_L=L^\ast/L$. It has the order $|\det L|$. The
group $A_L$ is equipped with the {\it discriminant bilinear form}
$b_L:A_L\times A_L\to \bq/\bz$ and the {\it discriminant quadratic
form} $q_L:A_L\to \bq/2\bz$ if $L$ is even. To get this forms, one
should extend the form of $L$ to the form on the dual lattice
$L^\ast$ with values in $\bq$.

\section{Isomorphisms between $M_X (v)$ for a general K3
surface $X$ and a primitive isotropic Mukai vectors $v$}
\label{secpicard1}

We consider algebraic K3 surfaces $X$ over $\bc$. We remind that a
non-singular projective  algebraic (or K\"ahlerian compact)
surface $X$ is a K3 surface if the canonical class $K_X$ is equal
to zero and the irregularity $q=\dim \Omega^1[X]=0$.

Further $N(X)$ denotes the Picard lattice of $X$, and $T(X)$ the
transcendental lattice of $X$.

We consider primitive isotropic Mukai vectors
\begin{equation}
v=(r,l,s),\ r\in \bn,\ s\in \bz,\  l\in N(X),\  l^2=2rs.
\label{mukaivector}
\end{equation}
on $X$. We denote by $Y=M_X(v)=M_X(r,l,s)$ the K3 surface which is the
minimal resolution of singularities of the moduli space of sheaves
over $X$ with the Mukai vector $v$. Clarify in Mukai \cite{Muk1}
--- \cite{Muk4} and Yoshioka \cite{Yoshioka1}.

In this section,  we call an {\it algebraic K3 surface}  to be {\it general} if
the Picard number $\rho (X)=\rk N(X)=1$ and the automorphsim group
of the transcendental periods of $X$ is trivial over $\bq$:
$Aut(T(X)\otimes \bq,\ H^{2,0}(X))=\pm 1$.

Here we consider the following question: When for a general algebraic
K3 surface $X$ and two its primitive isotropic Mukai vectors
$v_1=(r_1,l_1,s_1)$ and $v_2=(r_2,l_2,s_2)$, the moduli spaces
$M_X(v_1)$ and $M_X(v_2)$ are isomorphic?

\medskip

We have the following three {\it universal isomorphisms} between moduli
of shea\-ves over a K3 surface :

Let $D\in N(X)$.
Then one has the natural isomorphism given by the tensor product
$$
T_D:M_X(r,l,s)\cong M_X(r,l+rD,s+r(D^2/2)+D\cdot l),\ \
\E\mapsto \E\otimes \Oc(D).
$$
Moreover, here Mukai vectors
$$
v=(r,l,s),\ \ T_D(v)=(r,l+rD,s+r(D^2/2)+D\cdot l)
$$
have the same general common divisor and the same square under Mukai pairing.
In particular, they are primitive and isotropic simultaneously.

Taking $D=kH$ where $H$ is a hyperplane section and $k>0$,
using the isomorphisms $T_D$, we can always replace $M_X(r,l,s)$ by
an isomorphic
$M_X(r,l^\prime,s^\prime)$ where $l^\prime$ is ample, and then
${l^\prime}^2>0$. Thus, in our problem, we can further assume that
$v=(r,l,s)$ where $r>0$ and $l$ is ample. Then $l^2=2rs>0$ and $s>0$ either.
Further we assume that.

\medskip

For $r,\,s\,>0$, one has the isomorphism which is called {\it reflection}
$$
\delta:M_X(r,l,s)\cong M_X(s,l,r).
$$
E. g. see \cite{Muk2} and \cite{Tyurin1}, \cite{Tyurin2}, \cite{Yoshioka2}.
Thus, using the reflection, we can further assume that $0<r\le s$.

\medskip

For integers $d_1,\,d_2>0$ such that
$(d_1,d_2)=(d_1,s)=(r,d_2)=1$, one has the isomorphism
$$
\nu(d_1,d_2):M_X(r,l,s)\cong M_X(rd_1^2,ld_1d_2,sd_2^2)
$$
and its inverse $\nu(d_1,d_2)^{-1}$. See \cite{Muk2}, \cite{Muk3},
\cite{Nik3}, \cite{Mad-Nik5}. Using the isomorphisms
$\nu(d_1,d_2)$, $\nu(d_1,d_2)^{-1}$ and the reflection $\delta$,
we can always assume that the primitive isotropic Mukai vector
$v=(r,l,s)$ satisfies:
\begin{equation}
v=(r,l,s)\ \text{has}\ 0<r\le s,\ l^2=2rs\ \text{and}\ l\in N(X)\
\text{is primitive and ample.}
\label{redMukvec}
\end{equation}
We call such a primitive isotropic Mukai vector as a {\it reduced primitive
isotropic Mukai vector} (for $\rho(X)=1$).

We have the following result.

\begin{theorem}
\label{theoremunivred1}
Let $X$ be a general algebraic K3 surface, i. e.
$N(X)=\bz H$ where $H$ is a primitive polarization of $X$ and
$\Aut(T(X)\otimes \bq,H^{2,0}(X))=\pm 1$.  Let $v=(r,H,s)$ and
$v^\prime=(r^\prime,H, s^\prime)$ are two
reduced primitive isotropic Mukai vectors on $X$ (see \eqref{redMukvec}),
i. e. $0<r\le s$ and $0<r^\prime\le s^\prime$.

Then $M_X(v)\cong M_X(v^\prime)$ if and only if $v=v^\prime$, i. e.
$r^\prime=r$, $l^\prime =l$.
\end{theorem}

It follows that the described above universal (i. e. valid for all
algebraic K3 surfaces, even general ones) isomorphisms $T_D$,
$\delta$ and $\nu(d_1,d_2)$ are sufficient to find isomorphic
moduli of sheaves with primitive isotropic Mukai vectors for a
general K3 surface.

\begin{corollary}
\label{corollaryuniv}
Let $X$ be a general algebraic K3 surface and $v$, $v^\prime$ are
primitive isotropic Mukai vectors on $X$.

Then $M_X(v)\cong M_X(v^\prime)$ if and only if there exists an isomorphism
between $M_X(v)$ and $M_X(v^\prime)$ which is a composition of the universal
isomorphisms $T_D$, $\delta$ and $\nu(d_1,d_2)$.
\end{corollary}

\begin{proof}
Considerations below are very similar to much more general
and difficult calculations in (\cite{Nik5}, Sect. 2.3).

We have
$$
N(X)=\bz H=\{x\in H^2(X,\bz)\ |\ x\cdot H^{2,0}(X)=0\},
$$
and the transcendental lattice of $X$ is
$$
T(X)=N(X)^\perp_{H^2(X,\bz)}.
$$
The lattices $N(X)$ and $T(X)$ are orthogonal complements to one another in
the unimodular lattice
$H^2(X,\bz)$. We have that $N(X)\oplus T(X)\subset H^2(X,\bz)$ is only a
sublattice of a finite index. Here and in what follows $\oplus$ denotes the
orthogonal sum. Since $H^2(X,\bz)$ is unimodular and
$N(X)=\bz H$ is its primitive sublattice, there exists
$u\in H^2(X,\bz)$ such that $u\cdot H=1$.

We denote by $N(X)^\ast=\bz H/(2rs)
\subset
N(X)\otimes \bq$ and $T(X)^\ast\subset T(X)\otimes \bq$ the dual lattices.
Then $H^2(X,\bz)\subset N(X)^\ast\oplus T(X)^\ast$, and
$$
u=H/(2rs)\oplus t^\ast(H),\ t^\ast(H) \in T(X)^\ast.
$$
The element
$$
t^\ast (H) \mod T(X)\in T(X)^\ast/T(X)\cong \bz/2rs\bz
$$
is defined canonically by the primitive
element $H\in H^2(X,\bz)$. We evidently have
$$
H^2(X,\bz)=[N(X),T(X),u=H/(2rs)+t^\ast(H)]
$$
where $[\ \cdot\ ]$ means ``generated by'' $\cdot$\ .
The element
$t^\ast(H)\mod T(X)$ evidently distinguishes between different
polarized K3 surfaces with Picard number one and the same
transcendental periods. More exactly, for
another polarized K3 surface $(X^\prime , H^\prime )$ and its transcendental
periods $(T(X^\prime),H^{2,0}(X^\prime))$, the periods of $X$ and $X^\prime$
are isomorphic (and then $X\cong X^\prime$ by the Global Torelli Theorem
\cite{PS}) if and only if there exists an isomorphism
$\phi:T(X)\cong T(X^\prime)$ of the transcendental lattices such that
$(\phi\otimes \bc )(H^{2,0}(X))=H^{2,0}(X^\prime)$ and
$(\phi\otimes \bq) (t^\ast(H))\mod T(X)=t^\ast(H^\prime)\mod T(X^\prime)$.

Thus, the calculation of the periods of $X$ in terms its transcendental
periods is contained in the following statement.

\begin{proposition}
\label{propositionperX}
Let $(X,H)$ be a polarized K3 surface with a primitive
polarization $H$ such that $H^2=2rs$. Assume that
$N(X)=\bz H$ (i. e. $\rho(X)=1$).

Then
$$
H^2(X,\bz)=[N(X)=\bz H, T(X), H/(2rs)+t^\ast (H) ]
$$
where $t^\ast(H)\in T(X)^\ast$. The element $t^\ast (H)  \mod T(X)$ is defined
uniquely.

Moreover, $H^{2,0}(X)\subset T(X)\otimes \bc$.

(In general, when $\rho(X)\ge 1$,
one should replace $T(X)$ by $H^\perp_{H^2(X,\bz)}$.)
\end{proposition}

Let $Y=M_X(r,H,s)$. Let us calculate periods of $Y$.

We denote by
$$
\widetilde{H} (X,\bz)=H^0(X,\bz)+H^2(X,\bz)+H^4(X,\bz)=U\oplus H^2(X,\bz)
$$
(it is the direct sum) the Mukai lattice of $X$.
Here $U=\bz e_1+\bz e_2$ is the hyperbolic plane
where canonically $\bz e_1=H^0(X,\bz)$ and $\bz e_2=H^4(X,\bz)$
with the Mukai pairing $e_1^2=e_2^2=0$ and $e_1\cdot e_2=-1$.
Here $H^2(X,\bz)$ is the cohomology lattice of $X$ with the intersection
pairing. Here $\oplus$ denotes the orthogonal sum of lattices.
We have
\begin{equation}
v=re_1+se_2+H.
\label{v}
\end{equation}

By Mukai \cite{Muk2}, we have
\begin{equation}
H^2(Y,\bz)=v^\perp/\bz v,
\label{Mukai1}
\end{equation}
and  $H^{2,0}(Y)=H^{2,0}(X)$ by the canonical projection. This
defines periods of $Y$ and the isomorphism class of the K3 surface $Y$ (by
Global Torelli Theorem \cite{PS}). Let us calculate periods of $Y$ similarly to
Proposition \ref{propositionperX}.

Any element $f$ of $\widetilde{H}(X,\bz)$ can be uniquely written as
$$
f=xe_1+ye_2+\alpha H/(2rs)+\beta t^\ast,\ x,y,\alpha \in \bz,\ t^\ast\in
T(X)^\ast .
$$
We have $f\cdot v=-sx-ry+\alpha$, and $f\in v^\perp$ if and only if
$-sx-ry+\alpha=0$, and then
$$
f=xe_1+ye_2+(sx+ry)(H/(2rs))+\beta t^\ast.
$$
By Proposition \ref{propositionperX}, $f\in \widetilde{H}(X,\bz)$
if and only if
$t^\ast=(sx+ry)t^\ast (H) \mod T(X)$. Since $T(X)\subset v^\perp$, we can write
$$
f=xe_1+ye_2+(sx+ry)\left(H/(2rs)+t^\ast (H) \right)\mod T(X),\ x,y\in \bz\,.
$$
We denote
$$
c=(r,s),\ a=r/c,\ b=s/c.
$$
Then $(a,b)=1$. We have $h=-ae_1+be_2\in v^\perp$ and
$h^2=2ab=2rs/c^2$. Moreover, $h\perp T(X)$ and then $h \perp H^{2,0}(X)$.
Thus,
\begin{equation}
\label{definitionofh}
h\mod \bz v=-ae_1+be_2\mod \bz v
\end{equation}
gives an element of the Picard lattice $N(Y)$. We have
$$
e_1=\frac{v-ch-H}{2r},\ e_2=\frac{v+ch-H}{2s}.
$$
It follows that
\begin{equation}
f=\frac{sx+ry}{2rs}v+\frac{c(-sx+ry)}{2rs} h+(sx+ry)t^\ast (H) \mod T(X),\
x,y\in \bz\,.
\label{eleH2(Y)}
\end{equation}
Here $f\mod \bz v$ gives all elements of $H^2(Y,\bz)$ and
$H^{2,0}(Y)=H^{2,0}(X)\subset T(X)\otimes \bc$.

It follows that $f\mod \bz v\in T(Y)$ (where $\bz v$ gives the kernel of
$v^\perp$ and $H^2(Y,\bz)=v^\perp/\bz v$)
if and only if $-sx+ry=0$, equivalently
$-bx+ay=0$, equivalently (since $(a,b)=1$) $x=az$, $y=bz$ where $z\in \bz$,
and then $(sx+ry)t^\ast (H)=z(sa+rb)t^\ast (H)=z\,2abc\, t^\ast (H)$ where
$z\in \bz$. It follows that
\begin{equation}
T(Y)=[T(X),2abc \,t^\ast (H)].
\label{generalT(Y)}
\end{equation}
Since $t^\ast (H)  \mod T(X)$ has the order $2rs=2abc^2$ in $T(X)^\ast/T(X)\cong
\bz/2rs\bz$, it follows that $[T(Y):T(X)]=c$
(this is the result of Mukai, \cite{Muk2}).

By \eqref{eleH2(Y)}, \eqref{generalT(Y)}, we have
$f\perp H^{2,0}(Y)=H^{2,0}(X)$, equivalently
$f\mod \bz v\in N(Y)$, if and only if
$$
f=\frac{sx+ry}{2rs}v+\frac{c(-sx+ry)}{2rs} h
$$
where $sx+ry\equiv 0\mod 2abc$. Thus,
$acx+bcy\equiv 0\mod 2abc$ and $ax+by\equiv 0\mod 2ab$. Since $(a,b)=1$,
it follows that $x=b\widetilde{x}$, $y=a\widetilde{y}$ where
$\widetilde{x},\ \widetilde{y}\in \bz$, and
$\widetilde{x}+\widetilde{y}\equiv 0\mod 2$. Thus,
$\widetilde{y}=-\widetilde{x}+2k$ where $k\in \bz$. It follows that
$$
f=\frac{k}{c}v+(-\widetilde{x}+k)h, \ \ \widetilde{x},\, k\in \bz.
$$
Thus, $h\mod \bz v$ generates the Picard lattice $N(Y)$, and $h\mod \bz v$
can be considered as the polarization of $Y$
(or $-h\mod \bz v$ which does not matter from the point of view of periods and
the isomorphism class of $Y$).

Let us calculate $t^\ast (h)  \in T(Y)^\ast$. Then in \eqref{eleH2(Y)} we
should take an element
$f$ with $c(-sx+ry)/(2rs)=1/(2ab)$. Thus, $-sx+ry=c$ or $-bx+ay=1$. Then
$$
t^\ast (h)=(sx+ry)t^\ast (H) \mod T(Y).
$$
By \eqref{generalT(Y)}, $T(Y)^\ast=[T(X),ct^\ast (H) ]$ and
$T(Y)^\ast/T(Y)\cong \bz/2ab\bz$.

Thus, $t^\ast (h) =(bx+ay)\,(ct^\ast (H) \mod [T(X),2ab\,(ct^\ast (H) )]$
is defined by $m\equiv bx+ay\mod 2ab$. Since $-bx+ay=1$, we
have $m\equiv  2ay-1\equiv -1\mod 2a$ and $m\equiv 2bx+1\equiv 1\mod 2b$.
This defines $m\mod 2ab$ uniquely. We call such $m\mod 2ab$ as Mukai
element (compare with \cite{Muk3}).
Thus, $m(a,b)\mod 2ab$ is called {\it Mukai
element} if
\begin{equation}
m(a,b)\equiv -1\mod 2a\ \ \text{and}\ \ m(a,b) \equiv 1 \mod 2b.
\label{Mukele}
\end{equation}
Thus, $t^\ast (h) =m(a,b)\,ct^\ast (H)  \mod [T(X),2abc\,t^\ast (H) ]$.

Thus, finally, we finished the calculation of
periods of $Y$ in terms of periods of $X$ (see Proposition
\ref{propositionperX}).

\begin{proposition}
\label{propositionperY}
Let $(X,H)$ be a polarized K3 surface with a primitive
polarization $H$ such $H^2=2rs$, $r,s>0$. Assume that
$N(X)=\bz H$ (i. e. $\rho(X)=1$). Let $Y=M_X(r,H,s)$.
Let $c=(r,s)$ and $a=r/c$, $b=s/c$.

Then $N(Y)=\bz h$ where $h^2=2ab$,
$T(Y)=[T(X),2abc\,t^\ast (H) ]$, $T(Y)^\ast=[T(X), ct^\ast (H) ]$
and $t^\ast (h) \mod T(Y)=m(a,b)ct^\ast (H) \mod T(Y)$ where
$m(a,b)\mod 2ab$ is the Mukai element:
$m(a,b)\equiv -1 \mod 2a$, $m(a,b)\equiv 1\mod 2b$.
Thus,
$$
H^2(Y,\bz)=[N(Y),T(Y),h/(2ab)+t^\ast (h) ]=
$$
$$
[\bz h,\ [T(X),\, 2abc\,t^\ast (H) ],\ h/(2ab)+m(a,b)ct^\ast (H) ].
$$
(In general, when $\rho(X)\ge 1$, one should replace $T(X)$ by
$H^\perp_{H^2(X,\bz)}$ and $T(Y)$ by $h^\perp_{H^2(Y,\bz)}$.)
\end{proposition}

Now let us prove Theorem \ref{theoremunivred1}. We need to recover
$r$ and $s$ from periods of $Y$. By Proposition \ref{propositionperY},
we have $N(Y)=\bz h$ where $h^2=2ab$. Thus, we recover $ab$. Since
$c^2=2rs/2ab$, we recover $c$.

We have $(T(X)\otimes \bq,H^{2,0}(X))\cong (T(Y)\otimes \bq, H^{2,0}(Y))$.
Since $X$ is general, there exists only one such isomorphism up
to multiplication by $\pm 1$. It follows that there exists only one
(up to multiplication by $\pm 1$)
embedding $T(X)\subset T(Y)$ of lattices which identifies $H^{2,0}(X)$
and $H^{2,0}(Y)$. By Proposition \ref{propositionperY}, then
$t^\ast (h) \mod T(Y)=\widetilde{m}(a,b) c t^\ast (H) \mod T(Y)$ where
$\widetilde{m}(a,b)\equiv \pm m(a,b)\mod 2ab$ and $m(a,b)$ is
the Mukai element. Assume $p^\alpha |ab$ and $p^{\alpha+1}$ does not
divide $ab$ where $p$ is prime and $\alpha>0$. Then
$\widetilde{m}(a,b)\equiv \pm 1\mod 2p^\alpha$. Evidently, here only one
sing $\pm 1$ is possible, and we denote by $a$ the product of all such
$p^\alpha$ having $\widetilde{m}(a,b)\equiv -1 \mod 2p^\alpha$, and by by
$b$ the product of all such remaining $p^\alpha$ having
$\widetilde{m}(a,b)\equiv 1\mod 2p^\alpha$. If $a>b$, we should change
$a$ and $b$ places. Thus, we recover $a$ and $b$ and the reduced primitive
Mukai vector $(r,H,s)=(ac,H,bc)$ such that periods of $M_X(r,H,s)$ are
isomorphic to the periods of $Y$.

This finishes the proof.
\end{proof}

\medskip

\begin{remark}
\label{remarknegH}
 The same Propositions \ref{propositionperX}, \ref{propositionperY}
and the proofs  are valid for an algebraic K3 surface $X$
and a primitive element $H\in N(X)$ with $H^2=2rs\not=0$,
if one replaces $T(X)$ by the orthogonal complement
$H^\perp_{H^2(X,\bz)}$.
\end{remark}

\medskip

As an example of an application of Theorem \ref{theoremunivred1}, let
us consider the case when $M_X(r,l,s)\cong X$. It is known
(e.g. see \cite{Tyurin2}) that for $l\in N(X)$ and $\pm l^2>0$, one
has the Tyurin isomorphism
\begin{equation}
\label{Tyurinisom}
Tyu=Tyu(\pm l)\,:\, M_X(\pm l^2/2,l,\pm 1)\cong X \,.
\end{equation}
Existence of such an isomorphism follows at once from Global Torelli Theorem
for K3 surfaces \cite{PS} using Propositions \ref{propositionperX},
\ref{propositionperY} and Remark \ref{remarknegH}.

Thus, for a general K3 surface $X$ and a reduced primitive isotropic Mukai vector
$v=(r,H,1)$ where $r=H^2/2$, we have $M_X(r,H,1)\cong X$. By
Theorem \ref{theoremunivred1}, we then obtain the following result where
we also use the well-known fact that $\Aut (T(X),H^{2,0}(X))=\pm 1$  if
$\rho (X)=1$ (see \eqref{modXwithaut} below); it is sufficient to consider
the automorphism group over $\bz$ for this result.

\begin{corollary}
\label{corollaryunivred1congX}
Let $X$ be an algebraic K3 surface with $\rho (X)=1$, i. e.
$N(X)=\bz H$ where $H$ is a primitive polarization of $X$.
Let $v=(r,H,s)$ be
a reduced primitive isotropic Mukai vector on $X$ (see \eqref{redMukvec}),
i. e. $0<r\le s$.

Then $M_X(v)\cong X$ if and only if $v=(1,H,H^2/2)$, i. e.
$r=1$, $s=H^2/2$.
\end{corollary}

\section{Isomorphisms between $M_X (v)$ and $X$ for a general K3
surface $X$ with $\rho (X)=2$}
\label{secpicard2}

Here we consider general K3 surfaces $X$ with $\rho (X)=\rk N(X)=2$.
Here a K3 surface $X$ will be called general if the group of automorphisms
of the transcendental periods is trivial: $\Aut (T(X),H^{2,0}(X))=\pm 1$.

For $\rho (X)\ge 2$, we don't know when $M_X(v_1)\cong M_X(v_2)$ for
primitive isotropic Mukai vectors $v_1$ and $v_2$ on $X$. But we still
have the universal isomorphisms $T_D$, $D\in N(X)$, the reflection $\delta$,
the isomorphism $\nu(d_1,d_2)$ and the Tyurin isomorphism $Tyu$
considered in Section \ref{secpicard1}.
They are {\it universal isomorphisms,} i. e. they are
defined for all K3 surfaces, even with Picard number one.

First, we review results of \cite{Nik5} and \cite{Mad-Nik5} where for
general K3 surfaces $X$ with $\rho (X)=2$ all primitive isotropic Mukai
vectors $v$ with $M_X(v)\cong X$ were found. In particular, we know
when $M_X(v_1)\cong M_X(v_2)$ in the case when both moduli are
isomorphic to $X$. The result is that
$M_X(v)\cong X$ if and only if there exists such an isomorphism which is
a composition of the universal isomorphisms $\delta$, $T_D$ and
$\nu(d_1,d_2)$ between moduli of sheaves over $X$ and the Tyurin isomorphism
between moduli of sheaves over $X$ and $X$ itself. More exactly, the
results are as follows.

Using universal isomorphisms $T_D$, we can assume that the primitive
isotropic Mukai vector is
$$
v=(r,H,s),\ \ r>0,\ s>0,\ H^2=2rs.
$$
(We can even assume that $H$ is ample.)
We are interested in the case when $Y=M_X(r,H,s)\cong X$.

We denote $c=(r,s)$ and $a=r/c$, $b=s/c$. Then $(a,b)=1$. Let $H$
is divisible by $d\in \bn$ where $\wH=H/d$ is primitive in $N(X)$.
Primitivity of $v=(r,H,s)$ means that $(r,d,s)=(c,d)=1$. Since
$\wH^2=2abc^2/d^2$ is even, we have $d^2|abc^2$. Since
$(a,b)=(c,d)=1$, it follows that $d=d_ad_b$ where $d_a=(d,a)$,
$d_b=(d,b)$, and we can introduce integers
$$
a_1=\frac{a}{d_a^2},\ \ b_1=\frac{b}{d_b^2}\,.
$$
Then we obtain that $\wH^2=2a_1b_1c^2$.

Let $\gamma=\gamma (\wH)$ is defined by $\wH\cdot N(X)=\gamma
\bz$, i.e. $H\cdot N(X)=\gamma d\bz$. Clearly, $\gamma|\wH^2=2a_1b_1c^2$.

We denote
\begin{equation}
\label{n(v)1} n(v)=(r,s,d\gamma)=(r,s,\gamma).
\end{equation}
By Mukai \cite{Muk2}, we have $T(X)\subset T(Y)$, and
\begin{equation}
n(v)=[T(Y):T(X)]
\label{Mukai2}
\end{equation}
where $T(X)$ and $T(Y)$ are transcendental
lattices of $X$ and $Y$. Thus,
\begin{equation}
\label{n(v)2} Y\cong X\ \implies \
n(v)=(r,s,d\gamma)=(c,d\gamma)=(c,\gamma)=1.
\end{equation}

Assuming that $Y\cong X$ and then $n(v)=1$, we have
$\gamma|2a_1b_1$, and we can introduce
\begin{equation}
\gamma_a=(\gamma,a_1),\ \gamma_b=(\gamma,b_1),\
\gamma_2=\frac{\gamma}{\gamma_a\gamma_b}.
\end{equation}
Clearly, $\gamma_2|2$.

In (\cite{Nik5}, Theorem 4.4) the following general theorem had
been obtained (see its important particular cases in
\cite{Mad-Nik1}, \cite{Mad-Nik2} and \cite{Nik4}). In the theorem,
we use notations $c$, $a$, $b$, $d$, $d_a$, $d_b$, $a_1$, $b_1$
introduced above. The same notations $\gamma$, $\gamma_a$,
$\gamma_b$ and $\gamma_2$ as above are used if one replaces $N(X)$ by a
2-dimensional primitive sublattice
$N\subset N(X)$, e. g. $\wH\cdot N=\gamma \bz$,
$\gamma>0$. We denote $\det{N}=-\gamma\delta$ and $\bz f(\wH)$
denotes the orthogonal complement to $\wH$ in $N$.

\begin{theorem}
\label{maintheorem} Let $X$ be a K3 surface and $H$ a polarization
of $X$ such that $H^2=2rs$ where $r,s\in \bn$. Assume that the
Mukai vector $(r,H,s)$ is primitive. Let $Y=M_X(r,H,s)$ be the K3
surface which is the moduli of sheaves over $X$ with the isotropic
Mukai vector $v=(r,H,s)$. Let $\wH=H/d$, $d\in\bn$, be the
corresponding primitive polarization.

We have $Y\cong X$ if there exists $\wth_1 \in N(X)$ such that
$\wH$ and $\wth_1$ belong to a 2-dimensional primitive sublattice
$N \subset N(X)$ such that $\wH\cdot N=\gamma \bz$, $\gamma>0$,
$(c,d\gamma)=1$, and  the element $\wth_1$ belongs to the
$a$-series or to the $b$-series described below:

$\wth_1$ belongs to the $a$-series if
\begin{equation}
\wth_1^2=\pm 2b_1c, \ \ \wH\cdot \wth_1\equiv 0\mod
\gamma(b_1/\gamma_b)c,\ \ f(\wH)\cdot \wth_1\equiv 0\mod \delta
b_1c \label{aseries1}
\end{equation}
(where $\gamma_b=(\gamma,b_1)$);

$\wth_1$ belongs to the $b$-series if
\begin{equation}
\wth_1^2=\pm 2a_1c, \ \ \wH\cdot \wth_1\equiv 0\mod
\gamma(a_1/\gamma_a)c,\ \ f(\wH)\cdot \wth_1\equiv 0\mod\delta
a_1c \label{bseries1}
\end{equation}
(where $\gamma_a=(\gamma,a_1)$).

These conditions are necessary to have $Y\cong X$ if $\rho (X)\le
2$ and $X$ is a general K3 surface with its Picard lattice, i. e
the automorphism group of the transcendental periods
$(T(X),H^{2,0}(X))$ is $\pm 1$.
\end{theorem}

In \cite{Mad-Nik5}, Theorem \ref{maintheorem} was geometrically interpreted as
follows.

\begin{theorem}
\label{maintheorem2} Let $X$ be a K3 surface and $H$ a
polarization of $X$ such that $H^2=2rs$ where $r,s\in \bn$. Assume
that the Mukai vector $(r,H,s)$ is primitive. Let $Y=M_X(r,H,s)$
be the K3 surface which is the moduli of sheaves over $X$ with the
isotropic Mukai vector $v=(r,H,s)$. Let $\wH=H/d$, $d\in \bn$, be
the corresponding primitive polarization.

Assume that there exists $\wth_1 \in N(X)$ such that $\wH$ and
$\wth_1$ belong to a 2-dimensional primitive sublattice $N \subset
N(X)$ such that $\wH\cdot N=\gamma \bz$, $\gamma>0$,
$(c,d\gamma)=1$, and the element $\wth_1$ belongs to the
$a$-series or to the $b$-series described below:

$\wth_1$ belongs to the $a$-series if
\begin{equation}
\wth_1^2=\pm 2b_1c, \ \ \wH\cdot \wth_1\equiv 0\mod
\gamma(b_1/\gamma_b)c,\ \ f(\wH)\cdot \wth_1\equiv 0\mod \delta
b_1c \label{aseries2}
\end{equation}
(where $\gamma_b=(\gamma,b_1)$);

$\wth_1$ belongs to the $b$-series if
\begin{equation}
\wth_1^2=\pm 2a_1c, \ \ \wH\cdot \wth_1\equiv 0\mod
\gamma(a_1/\gamma_a)c,\ \ f(\wH)\cdot \wth_1\equiv 0\mod\delta
a_1c \label{bseries2}
\end{equation}
(where $\gamma_a=(\gamma,a_1)$).

Then we have:

If $\wth_1$ belongs to the $a$-series, then
\begin{equation}
\wth_1=d_2\wH+b_1c D\ \text{for some\ } d_2\in \bn,\  D\in N,
\label{aseriesnew}
\end{equation}
which defines the isomorphism
\begin{equation}
Tyu(\pm \wth_1)\cdot T_D\cdot \nu(1,d_2)\cdot \delta\cdot
\nu(d_a,d_b)^{-1}: Y=M_X(r,H,s)\cong X . \label{aseriesisomnew}
\end{equation}

If $\wth_1$ belongs to the $b$-series, then
\begin{equation}
\wth_1=d_2\wH+a_1cD\ \text{for some\ } d_2\in \bn,\ D\in N,
\label{bseriesnew}
\end{equation}
which defines the isomorphism
\begin{equation}
Tyu(\pm \wth_1)\cdot T_D\cdot \nu(1,d_2)\cdot \nu(d_a,d_b)^{-1}:
Y=M_X(r,H,s)\cong X . \label{bseriesisomnew}
\end{equation}
\end{theorem}

Since conditions of Theorems \ref{maintheorem}, \ref{maintheorem2}
are necessary for general K3 surfaces
with $\rho (X)\le 2$, we obtain

\begin{theorem}
\label{maintheorem3}
Let $X$ be a K3 surface with a polarization
$H$ such that $H^2=2rs$, $r,s\ge 1$, the Mukai vector $(r,H,s)$ be
primitive, and $Y=M_X(r,H,s)$ be the moduli of sheaves over $X$
with the isotropic  Mukai vector $(r,H,s)$. Assume that
$\rho(X)\le 2$ and $X$ is general with its Picard lattice (i. e.
the automorphism group of the transcendental periods $Aut(T(X),
H^{2,0}(X))=\pm 1$). Let $\wH=H/d$, $d\in \bn$, be the corresponding
primitive polarization.

Then $Y=M_X(r,H,s)$ is isomorphic to $X$ if and only if there
exists $d_2\in \bn$ and $D\in N=N(X)$ such that

either
\begin{equation}
\wth_1=d_2\wH+b_1c D\ \text{has}\ \wth_1^2=\pm 2b_1c,\
\label{aseriesnew2}
\end{equation}
which defines the isomorphism
\begin{equation}
Tyu(\pm \wth_1)\cdot T_D\cdot \nu(1,d_2)\cdot \delta\cdot
\nu(d_a,d_b)^{-1}: Y=M_X(r,H,s)\cong X,\label{aseriesisomnew2}
\end{equation}

or
\begin{equation}
\wth_1=d_2\wH+a_1cD\ \text{has}\ \wth_1^2=\pm 2a_1c,
\label{bseriesnew2}
\end{equation}
which defines the isomorphism
\begin{equation}
Tyu(\pm \wth_1)\cdot T_D\cdot \nu(1,d_2)\cdot \nu(d_a,d_b)^{-1}:
Y=M_X(r,H,s)\cong X . \label{bseriesisomnew2}
\end{equation}
\end{theorem}

Theorem \ref{theoremunivred1} clarifies appearance of the
isomorphisms $T_D$, $\delta$, $\nu(d_1,d_2)$ and $Tyu$ in these
results for Picard number 2. They are universal and exist for all
K3 surfaces; moreover, they are all isomorphisms which one needs  to
distinguish isomorphic moduli $M_X(v)$ for isotropic Mukai vectors
$v$ on a general K3 surface $X$. Thus, appearance of the
isomorphisms $T_D$, $\delta$, $\nu(d_1,d_2)$ and $Tyu$ is very
natural in the results above.

\begin{remark}
\label{remarkrho=1}
For The Picard number $\rho(X)=1$,
Theorems \ref{maintheorem}, \ref{maintheorem2} and \ref{maintheorem3} are
formally equivalent to Corollary \ref{corollaryunivred1congX}. Really,
for $\rho(X)=1$ we have $\gamma=2a_1b_1c^2$. Thus, $(\gamma,c)=1$ implies
that $c=1$. Then $\gamma=2a_1b_1$ and  $\gamma_2=2$,
$\gamma_a=a_1$, $\gamma_b=b_1$. Conditions of Theorem \ref{maintheorem}
can be satisfied only for $\wth_1=\wH$ which implies that $a_1=1$ for
the $a$-series, and $b_1=1$ for the $b$-series
(one can formally put $f(\wH)=0$).

Thus, for $\rho (X)=1$ and general $X$,
we have $Y\cong X$ if and only if $c=1$ and
either $a_1=1$ or $b_1=1$. This is equivalent to Corollary
\ref{corollaryunivred1congX}.
\end{remark}

\medskip

Under conditions of Theorem \ref{maintheorem}, let
us assume that for a primitive 2-dimensional sublattice
$N\subset N(X)$ an element $\wth_1\in N$ with $\wth_1^2=\pm 2b_1c$
belongs to the $a$-series. This is equivalent to the condition
\eqref{aseriesnew} of Theorem \ref{maintheorem2}.
Replacing $\wth_1$ by $-\wth_1$ if necessary, we see that
\eqref{aseriesnew} is equivalent to
\begin{equation}
\label{aseriesnew1}
\wth_1=d_2\wH+b_1c D,\ \ \ d_2\in \bz,\ D\in N.
\end{equation}
Since $\wH$ is primitive, the lattice $N$ has a basis $\wH$, $D\in N$, i. e.
$N=[\wH,\,D]$. Since $\wH\cdot N=\gamma \bz$ where $(\gamma,c)=1$, the matrix
of $N$ in this basis is
\begin{equation}
\left(
\begin{array}{cc}
\wH^2         & \wH\cdot D\\
\wH\cdot D    & D^2
\end{array}
\right)
\ =\
\left(
\begin{array}{cc}
2a_1b_1c^2   & \gamma k\\
\gamma k    & 2t
\end{array}
\right)
\label{matrixofN}
\end{equation}
where $k,t\in \bz$ and $\gamma|2a_1b_1$, $(\gamma,c)=1$ and
$(2a_1b_1c^2/\gamma,k)=1$.

The condition of $a$-series \eqref{aseriesnew1}
is then equivalent to existence of
$\wth_1\in [\wH,b_1c N]=[\wH,b_1c D]$ with $\wth_1^2=\pm 2b_1c$.
Thus, the lattice $N_1=[\wH,b_1cD]$ with the matrix
\begin{equation}
\left(
\begin{array}{cc}
2a_1b_1c^2   & b_1c\gamma k\\
b_1c\gamma k    & b_1^2c^2 2t
\end{array}
\right)
\label{matrixofN1}
\end{equation}
must have $\wth_1$ with $\wth_1^2=\pm 2b_1c$.
Writing $\wth_1$ as $\wth_1=x\wH+yb_1cD$, we obtain
that the quadratic equation
$a_1cx^2+\gamma kxy+b_1cty^2=\pm 1$ must have an integral solution.
Similarly, for $b$-series we obtain the equation
$b_1cx^2+\gamma kxy+a_1cty^2=\pm 1$. Thus, we finally obtain a very
elementary reformulation of the results above.

\begin{lemma}
\label{lemmaequations}
For the matrix \eqref{matrixofN} of the lattice
$N$ in Theorems \ref{maintheorem}, \ref{maintheorem2} and
\ref{maintheorem3}, the conditions of $a$-series are equivalent to
existence of an integral solution of the equation
\begin{equation}
a_1cx^2+\gamma kxy+b_1cty^2=\pm 1,
\label{aseriesequation}
\end{equation}
and for $b$-series of the equation
\begin{equation}
b_1cx^2+\gamma kxy+a_1cty^2=\pm 1.
\label{bseriesequation}
\end{equation}
\end{lemma}

This calculation has a very important corollary. Let us assume that
a prime $p|\gamma_b=(\gamma,b_1)$. Then for the equation
\eqref{aseriesequation} we obtain a congruence
$a_1cx^2\equiv \pm 1\mod p$. Thus, $\pm a_1c$ is a square $\mod p$.
Similarly, for the equation \eqref{bseriesequation}, we obtain that
$\pm b_1c$ is a square $\mod p$ for a prime $p|\gamma_a=(\gamma,a_1)$.

Thus, we obtain an important necessary condition of $Y=M_X(v)\cong X$
for $\rho (X)=2$.

\begin{theorem}
\label{theoremnecessary}
Let $X$ be a K3 surface with a polarization
$H$ such that $H^2=2rs$, $r,s\ge 1$, the Mukai vector $(r,H,s)$ be
primitive, and $Y=M_X(r,H,s)$ be the moduli of sheaves over $X$
with the isotropic  Mukai vector $(r,H,s)$. Assume that
$\rho(X)\le 2$ and $X$ is general with its Picard lattice (i. e.
the automorphism group of the transcendental periods $Aut(T(X),
H^{2,0}(X))=\pm 1$). Let $\wH=H/d$, $d\in \bn$, be the corresponding
primitive polarization, $\wH\cdot N(X)=\gamma \bz$ and $(\gamma,c)=1$.

Then $Y=M_X(r,H,s)\cong X$ implies that for one of $\pm$
either
\begin{equation}
\forall\  p|\gamma_b \implies \left(\frac{\pm a_1c}{p}\right)=1
\label{aseriesnescond}
\end{equation}
or
\begin{equation}
\forall p|\gamma_a \implies \left(\frac{\pm b_1c}{p}\right)=1.
\label{bseriesnescond}
\end{equation}
Here $p$ means any prime, and $\left(\frac{x}{2}\right)=1$ means
that $x\equiv 1 \mod 8$.

Thus, if
\begin{equation}
\forall\ \pm\
\left(\exists\  p|\gamma_b\ \text{such that}\
\left(\frac{\pm a_1c}{p}\right)=-1\ \text{and}\
\exists p|\gamma_a \text{such that} \left(\frac{\pm b_1c}{p}\right)=-1\right)
\label{seriesnescondopp}
\end{equation}
then $Y=M_X(r,H,s)$ is not isomorphic
to $X$ for a general (for its Picard lattice) K3 surface $X$ with
$\rho (X)\le 2$.
\end{theorem}

\begin{example}
\label{example5,13}
Assume that $a_1=5$, $b_1=13$, $c=1$ and $\gamma=5\cdot 13$
(or $\gamma=2\cdot 5\cdot 13$). Obviously, then \eqref{seriesnescondopp}
is valid.  Thus, for
$$
v=(5, H,13),\ \ H^2=2\cdot 5\cdot 13,\ \text{and}\
\gamma=5\cdot 13\ \text{or}\
2\cdot 5\cdot 13
$$
(then $H$ is always primitive),
for any general K3 surface $X$ with $\rho(X)=2$ and any $H\in N(X)$
with $H^2=2\cdot 5\cdot 13$ and $H\cdot N(X)=\gamma \bz $,
the moduli $Y=M_X(v)$ are not isomorphic to $X$.

There are a lot of such Picard lattices given by \eqref{matrixofN}.
\end{example}

In \cite{Nik4}, it is shown that any primitive isotropic Mukai
vector $v=(r,H,s)$ with $H^2=2rs$ and $\gamma=1$ can be realized by
a general K3 surface with Picard number 2. It is possible that
Theorem \ref{theoremnecessary} gives all necessary conditions to
have similar result for any $\gamma$. We hope to consider this problem
later.

\medskip

Importance of these results for general K3 surfaces $X$ with $\rho
(X)=2$ is that these results describe {\it all divisorial
conditions on moduli of polarized K3 surfaces which imply that
$Y=M_X(r,H,s)\cong X$.} Let us consider the corresponding simple
general arguments.

It is well-known (see \cite{Nik1} and \cite{Nik2.5} where, it
seems, it was observed first) that $Aut(T(X),H^{2,0}(X))\cong C_m$
is a finite cyclic group of the order $m$, and its representation
in $T(X)\otimes \bq$ is the sum of irreducible representations of
the dimension $\phi (m)$ ($\phi (m)$ is the Euler function), and
$H^{2,0}(X)$ is a line in one of eigen-spaces of $C_m$. In
particular, $\phi (m)|rk T(X)$ and the dimension of moduli of
these $X$ is equal to
\begin{equation}
\label{modXwithaut}
\dim Mod(X)=\rk T(X)/\phi (m)-1
\end{equation}
if $m>2$. If $m=1,2$, then $\dim Mod(X)=\rk T(X)-2$.

Let us consider polarized K3 surfaces $(X,H)$ with the
polarization $H^2=2rs$ and a primitive Mukai vector
$(r,H,s)$, $r,s>0$. Let us assume that $Y=M_X(r,H,s)\cong X$.

If $\rho (X)=1$, then $\rk T(X)=21$ and $\phi (m)|21$. It follows
that $m=1$ or $m=2$ because $21$ is odd. Thus
$Aut(T(X),H^{2,0}(X))=\pm 1$, and then $c=1$ and either $a_1=1$ or
$b_1=1$ by Corollary \ref{corollaryunivred1congX} (or Remark
\ref{remarkrho=1}). By specialization principle (see Lemma 2.1.1
in \cite{Nik5}), then $Y\cong M_X(r,H,s)$ for all K3 surfaces $X$
and the Mukai vector with these invariants:
\begin{equation}
c=1,\  \text{and either\ } a_1=1\ \text{or}\  b_1=1.
\label{conditionrho=1}
\end{equation}

Now let us assume that $(r,H,s)$ does not satisfy
\eqref{conditionrho=1}, but $Y=M_X(r,H,s)$ $\cong X$. By Corollary
\ref{corollaryunivred1congX}, then $\rho (X)\not=1$. Then $\rho (X)\ge 2$
and $\dim Mod(X)\le 20-\rho (X)\le 18$. Thus, a divisorial condition on
moduli or polarized K3 surfaces $(X,H)$
to have $Y=M_X(r,H,s)\cong X$  means that $\rho (X)=2$
for a general K3 surface satisfying this condition.
All these conditions are described by the isomorphism
classes of $H\in N(X)$ where $\rk N(X)=2$ and $H\in N(X)$ satisfies
the equivalent (for this case) Theorems \ref{maintheorem}, \ref{maintheorem2}
or \ref{maintheorem3}.
If $H\in N\subset N(X)$ is a primitive sublattice of the rank two
and $H\in N$ satisfies equivalent Theorems \ref{maintheorem}
\ref{maintheorem2}, then $Y=M_X(r,H,s)\cong X$ by the specialization
principle. This means that $X$ belongs to the closure of
the divisor defined by the moduli of polarized K3 surfaces $(X^\prime,H)$
with the Picard lattice $N(X^\prime)=N$ of the rank two. Thus,
$Y^\prime=M_X(r,H,s)\cong X^\prime$ because $X^\prime $ satisfies
the divisorial condition $H\in N$ where $H\in N\subset N(X^\prime)$.

By Theorem \ref{theoremnecessary} we obtain

\begin{theorem}
\label{divconditionsnes} Let for $r,s\ge 1$,
$$
v=(r,H,s),\ H^2=2rs
$$
be a type of a primitive isotropic Mukai vector, and
$\gamma|2a_1b_1$ and $(\gamma,c)=1$.

Then, if \eqref{seriesnescondopp} is valid,
there does not exist a divisorial
condition on moduli of polarized K3 surfaces $(X,H)$ which implies
$Y=M_X(r,H,s)\cong X$ and $H\cdot N(X)=\gamma \bz$. Thus, these K3
surfaces have codimension at least 2 in 19-dimensional
moduli of polarized K3 surfaces $(X,H)$.

For example, this is valid for $r=5$, $s=13$ (then $H$ is primitive and
$d=1$), and $\gamma=5\cdot 13$ (or $\gamma=2\cdot 5\cdot 13$).
\end{theorem}

In the section below, we will show that the numerical example of
Theorem \ref{divconditionsnes} can be satisfied by K3 surfaces
$X$ with $\rho (X)=3$. Thus, this K3 surfaces define a 17-dimensional
submanifold in the moduli of polarized K3 surfaces. This submanifold
cannot be extended to a divisor in moduli
preserving the condition $Y=M_X(r,H,s)\cong X$.

\section{Isomorphisms between $M_X (v)$ and $X$ for a general K3
surface $X$ with $\rho (X)\ge 3$}
\label{secpicard3}

Here we show that it is interesting and non-trivial to
generalize results of the previous section to $\rho (X)\ge 3$.

Let $K=[e_1,e_2, (e_1+e_2)/2]$ be a negative definite 2-dimensional
lattice with $e_1^2=-6$, $e_2^2=-34$ and $e_1\cdot e_2=0$. Then
$\left((e_1+e_2)/2\right)^2=(-6-34)/4=-10$ is even, and the lattice
$K$ is even. Since $6x^2+34y^2=8$ has no integral solutions, it follows that
$K$ has no elements $\delta \in K$ with $\delta^2=-2$. Let us consider the
lattice
$$
S=\bz H\oplus K
$$
which is the orthogonal sum of $\bz H$ with $H^2=2\cdot 5\cdot 13$
and the lattice $K$. By standard results about K3 surfaces, there exists
a polarized K3 surface $(X,H)$ with the Picard lattice $S$ and
the polarization $H\in S$. We then have
$H\cdot S=2\cdot 5\cdot 13\,\bz$. Thus, $\gamma=2\cdot 5\cdot13$.

Let $Y=M_X(5,H,13)$. We have the following result. Perhaps, it gives the main
result of the paper.

\begin{theorem}
\label{theoremrho=3}
For any polarized K3 surface $(X,H)$ with $N(X)=S$ where $S$ is the hyperbolic
lattice of the rank 3 defined above,
one has $Y=M_X(5,H,13)\cong X$ which gives
a 17-dimensional moduli $M_S$ of polarized K3 surfaces
$(X,H)$ with $Y=M_X(5,H,13)\cong X$.

On the other hand, $M_S$ is not contained in any 18-dimensional moduli
$M_N$ of polarized K3 surfaces $(X^\prime,H)$ where
$H\in N(X^\prime)=N\subset S$, $\rk N=2$ and
$M_{X^\prime}(5,H,13)\cong X^\prime$. Thus, $M_S$ is not
defined by any divisorial condition on moduli of polarized K3 surfaces
$(X,H)$ (or, it is not a specialization of) implying $M_X(5,H,13)\cong X$.
\end{theorem}

\begin{proof} For this case, $c=(5,13)=1$ and $(\gamma,c)=1$. By Mukai
results \eqref{Mukai1} and \eqref{Mukai2}, then the transcendental periods
$(T(X),H^{2,0}(X))$ and $(T(Y),H^{2,0}(Y))$
are isomorphic. The discriminant group $A_S=S^\ast/S$ of the lattice
$S=T(X)^\perp$ is a cyclic group
$\bz/(2\cdot 5\cdot 13\cdot 3\cdot 17)$. Thus, the minimal number $
l(A_S)$ of generators of $A_S$ is one.
Thus, $l(A_S)\le \rk S-2$. By Theorem 1.14.4 in \cite{Nik2}, a
primitive embedding of $T(X)$ into the cohomology lattice of K3 (which
is an even unimodular lattice of signature $(3,19)$) is then unique, up
to isomorphisms.
It follows that the isomorphism between transcendental periods of
$X$ and $Y$ can be extended to an isomorphism of periods of $X$ and $Y$.
By Global Torelli Theorem for K3 surfaces \cite{PS}, the K3 surfaces
$X$ and $Y$ are isomorphic. (These considerations are now standard.)

Let $H\in N\subset S$ be a primitive sublattice with $\rk N=2$.
Since $H\cdot S=H\cdot H\bz =2\cdot 5\cdot 13\bz$, it follows that
$H\cdot N=2\cdot 5\cdot 13\bz$, and the invariant $\gamma=2\cdot
5\cdot 13$ is the same for any sublattice $N\subset S$ containing
$H$. By Theorem \ref{divconditionsnes}, then $M_{X^\prime}
(r,H,s)$ is not isomorphic to $X^\prime$ for any general K3
surface $(X^\prime,H)$ with $N(X^\prime)=N$.

This finishes the proof.
\end{proof}

Similar arguments can be used to prove the following general statement
for $\rho (X)\ge 12$ which shows that there are many cases when
$Y=M_X(r,H,s)\cong X$ which don't follow from divisorial conditions on
moduli. Its first statement is well-known (e. g. see
Proposition 2.2.1 in \cite{Mad-Nik1}).

\begin{theorem}
\label{theoremrhomore12} Let $(X,H)$ be a polarized K3 surface with
$\rho (X)\ge 12$, and for $r,s\ge 1$ let
$(r,H,s)$ be a primitive isotropic Mukai vector on $X$, i. e.
$H^2=2rs$ and $(c,d)=1$. Assume that $H\cdot N(X)=\gamma \bz$.

Then $Y=M_X(r,H,s)\cong X$ if $(\gamma,c)=1$ (Mukai necessary condition).

On the other hand, if \eqref{seriesnescondopp} satisfies, the
isomorphism $Y=M_X(r,H,s)\cong X$ does not follow from any
divisorial condition on moduli of polarized K3 surfaces. I. e. for
any primitive 2-dimensional sublattice $H\in N\subset N(X)$, there
exists a polarized K3 surface $(X^\prime,H)$ with $N(X^\prime)=N$ such that
$Y^\prime=M_{X^\prime}(r,H,s)$ is not isomorphic to $X^\prime$.
\end{theorem}

\begin{proof} Since $\rho (X)\ge 12$, then $\rk T(X)\le 22-12=10$ and
$l(A_{T(X)})\le \rk T(X)=10$. Since $N(X)$ and $T(X)$ are orthogonal
complements to one another in the unimodular lattice $H^2(X,\bz)$,
it follows that $A_{N(X)}\cong A_{T(X)}$ and
$l(A_{N(X)})\le 10\le \rk N(X)-2$. By Theorem 1.14.4 in \cite{Nik2},
a primitive embedding of $T(X)$ into the cohomology lattice of K3 is then
unique up to isomorphisms. Like in the proof of Theorem \ref{theoremrho=3},
it follows that $Y\cong X$.

Let us prove the second statement. Since $H\cdot N(X)=\gamma \bz$ and
$H\in N\subset N(X)$, it follows that $H\cdot N(X^\prime)=\gamma (N)\bz$
where $\gamma|\gamma(N)$. If $(c,\gamma (N))>1$, then $Y^\prime$
is not isomorphic to $X$ because
$[T(Y^\prime):T(X^\prime)]=(c,\gamma (N))>1$
by Mukai's result \eqref{Mukai2}. Assume $(c,\gamma (N))=1$. Obviously,
\eqref{seriesnescondopp} for $\gamma$ implies \eqref{seriesnescondopp} for
$\gamma (N)$. Let $X^\prime$ be a general K3 surface with $N(X^\prime )=N$.
By Theorem  \ref{theoremnecessary}, then $Y^\prime=M_{X^\prime}(r,H,s)$
is not isomorphic to $X^\prime$.

This finishes the proof.
\end{proof}

Theorems \ref{theoremrho=3} and \ref{theoremrhomore12} can be unified to
the following the most general (known) statement: when
$Y=M_X(r,H,s)\cong X$ for any primitive isotropic Mukai vector on $X$
satisfying Mukai's necessary condition.

\begin{theorem}
\label{theoremuniqtrper} Let $X$ be a K3 surface, the Picard
lattice $N(X)$ is unique in its genus, and the natural
homomorphism
$$
O(N(X))\to O(q_{N(X)})
$$
is surjective where $q_{N(X)}$ is the discriminant quadratic form
of $N(X)$.  Equivalently, any isomorphism of the transcendental
periods of $X$ and another K3 surface can be extended to the
isomorphisms of periods of $X$ and the other K3 surface.

Then for any primitive isotropic Mukai vector $v=(r,H,s)$ on $X$ such that
$(c,\gamma)=1$ (Mukai necessary condition), one has $Y=M_X(r,H,s)\cong X$.

On the other hand, if $X$ is general, i. e. $\Aut(T(X),H^{2,0}(X))=\pm 1$,
and \eqref{seriesnescondopp} satisfies,
then the isomorphism $Y=M_X(r,H,s)\cong X$ does not follow from any
divisorial condition on moduli of polarized K3 surfaces $(X,H)$. I. e. for
any primitive 2-dimensional sublattice $H\in N\subset N(X)$, there
exists a polarized K3 surface $(X^\prime,H)$ with $N(X^\prime)=N$ such that
$Y^\prime=M_{X^\prime}(r,H,s)$ is not isomorphic to $X^\prime$.
\end{theorem}

The results of Section \ref{secpicard2} and these results suggest the
following general concepts.

Let $r\in \bn $ and $s\in \bz$. We formally put $H^2=2rs$ and introduce
$c=(r,s)$ and $a=r/c$, $b=s/c$. Let $d\in \bn$, $(d,c)=1$ and $d^2|ab$.
We call
\begin{equation}
(r,H,s),\ H^2=2rs,\ d,
\label{typemukaivector}
\end{equation}
the {\it type of a primitive isotropic Mukai vector of K3.} Clearly,
{\it a Mukai vector of
the type \eqref{typemukaivector} on a K3 surface $X$} is just an element
$H\in N(X)$ such that $H^2=2rs$ and $\wH=H/d$ is primitive. Like above,
we introduce $d_a=(d,a)$, $d_b=(d,b)$ and put $a_1=a/d_a^2$, $b_1=b/d_b^2$.
Then $\wH^2=2a_1b_1c^2$.

Let $N$ be an even lattice which can be primitively embedded into a Picard
lattice of some algebraic K3 surface (equivalently, there exists a K\"ahlerian
K3 surface with this Picard lattice).
This is equivalent for $N$ to be either
negative definite, or semi-negative definite with 1-dimensional kernel,
or hyperbolic (i. e. $N$ has the signature $(1, \rho-1)$),
and to have a primitive embedding into an even unimodular lattice of
the signature $(3,19)$.
Further, we call $N$ as an {\it abstract K3 Picard lattice} (or just a
K3 Picard lattice).
Let $H\in N$. We call $H\in N$ as a {\it polarized (abstract) K3 Picard
lattice,} in spite of $H^2$ can be non-positive.
We consider such pairs up to natural isomorphisms. Another
polarized K3 Picard lattice $H^\prime\in N^\prime$ is called
isomorphic to $H\in N$ if there exists an isomorphism $f:N\cong N^\prime$
of lattices such that $f(H)=H^\prime$.

\begin{definition}
\label{defcritical} Let us fix a type \eqref{typemukaivector} of a
primitive isotropic Mukai vector of K3. A polarized K3 Picard
lattice $H\in N$ is called {\bf critical for correspondences of a
K3 surface with itself via moduli of sheaves for the type
\eqref{typemukaivector} of Mukai vector} (further we
abbreviate this as {\bf $H\in N$ is a critical polarized K3 Picard
lattice for the type \eqref{typemukaivector}}) if $H^2=2rs$ and
$\wH=H/d\in N$ is primitive; moreover the conditions (a) and (b)
below satisfy:

(a) for any K3 surface $X$ such that $H\in N\subset N(X)$ is a primitive
sublattice, one has
$Y=M_X(r,H,s)\cong X$.

(b) the condition (a) above is not valid if one replaces $H\in N$
by $H\in N_1$ for any primitive sublattice $H\in N_1\subset N$ of $N$
of strictly smaller rank $\rk N_1<\rk N$.
\end{definition}

\medskip

On the one hand, in (\cite{Nik5}, Theorem 2.3.3), for a polarized
K3 Picard lattice $H\in N$, the criterion is given for a general
(and then any) K3 surface with $H\in N= N(X)$ to have
$Y=M_X(r,H,s)\cong X$. On the other hand, by the specialization
principle (Lemma 2.1.1 in \cite{Nik5}), if this criterion is
satisfied, then $Y=M_X(r,H,s)\cong X$ for any K3 surface $X$ such
that $H\in N\subset N(X)$ is a primitive sublattice. Thus, for the
problem of describing in terms of Picard lattices, of all K3
surfaces $X$ such that $Y=M_X(r,H,s)\cong X$, the main problem is
as follows.

\begin{problem}
\label{problemcritical}
For a given type of a primitive isotropic Mukai vector
\eqref{typemukaivector} of K3, describe all {\bf critical polarized
K3 Picard lattices $H\in N$} (for the problem of correspondences of a K3
surface with itself via moduli of sheaves).
\end{problem}

Now we have the following examples of solution of this problem.

\medskip

By \eqref{Tyurinisom}, or Corollary \ref{corollaryunivred1congX},
or Remark \ref{remarkrho=1}, we have classification of
critical polarized K3 Picard lattices of the rank one.

\begin{example}
\label{examplerho=1}
For the type $(r,H,s)$, $H^2=2rs$, $d$
where $c=1$ and either $a_1=1$ or $b_1=\pm 1$, we obtain that
$N=\bz \wH$ where $\wH^2=2a_1b_1$ gives all critical polarized K3
Picard lattices $H=d\wH \in N$ of the rank one.
\end{example}

\begin{example}
\label{examplerho=2}
For the type of Mukai vector which is different from Example
\ref{examplerho=1}, classification of the critical polarized
K3 Picard lattices of the rank 2 is given by equivalent
Theorems \ref{maintheorem},
\ref{maintheorem2} or \ref{maintheorem3}.
\end{example}

\begin{example}
\label{examplerho=3} For the Mukai vector of the type $(5,H,13)$ with
$H^2=2\cdot 5\cdot 13$ and $d=1$, the polarized Picard lattice $H\in S$
of Theorem \ref{theoremrho=3} is critical of the rank $\rk S$=3, by
Theorem \ref{theoremrho=3}. Obviously, there are plenty of similar examples.
It would be very interesting and non-trivial
to find all critical polarized K3 Picard lattices $H\in S$ of the rank 3.
\end{example}

\begin{example}
\label{examplerhomore3} By Theorem \ref{theoremrhomore12}, we should
expect that there exist critical polarized K3 Picard lattices of the
rank more than 3. On the other hand, the same Theorem \ref{theoremrhomore12}
gives that the rank of a critical polarized K3 Picard lattice is less or
equal to 12.
\end{example}

We have

\begin{theorem}
\label{theoremrankcrit} For any type $(r,H,s)$, $H^2=2rs$ and $d$
of a primitive isotropic Mukai vector of K3, the rank
of a critical polarized K3 Picard lattice $H\in N$ is not more than 12:
we have $\rk N\le 12$.
\end{theorem}

\begin{proof} Let $H\in N$ be a critical polarized K3 Picard lattice of
this type and $\rk N\ge 13$. Let us take any primitive
sublattice $H\in N^\prime\subset N$ of the $\rk N^\prime=12$ such that
$\wH\cdot N^\prime =\wH\cdot N$. Obviously, it does exist.
Let $X$ be an algebraic K3 surface such that $H\in N^\prime\subset N(X)$.
Then $\rk N(X)\ge 12$ and $Y=M_X(r,H,s)\cong X$ by Theorem
\ref{theoremrhomore12}.

Then the condition (b) of Definition \ref{defcritical} is not satisfied,
and we get a contradiction. Thus, $\rk N\le 12$.

This finishes the proof.
\end{proof}

It would be very interesting to give an exact estimate for the rank of
critical polarized K3 Picard lattices.

\begin{problem}
\label{problemcritrank}
For a given type \eqref{typemukaivector} of a primitive isotropic Mukai
vector of K3, give the exact estimate of the rank $\rk N$ of
a critical polarized K3 Picard lattices $H\in N$ of this type (for the
problem of isomorphisms of K3 surfaces with itself).
\end{problem}

Now we don't know the answer to this problem for any type
\eqref{typemukaivector} different from Example \ref{examplerho=1}.

\section{Compositions of correspondences of a K3 surface
with itself via moduli of sheaves. General Problem of classification of
correspondences of K3 surface with itself via moduli of sheaves}
\label{seccompositions}

Here we want to interpret the results above in terms of the action of
correspondences on 2-dimensional cohomology lattice of a K3 surface.
Moreover, we try to formulate a general problem of classification of
correspondences of a K3 surface with itself via moduli of sheaves.

\medskip

Let $v=(r,H,s)$ be a primitive isotropic Mukai vector on a K3 surface $X$
and $Y=M_X(r,H,s)$. We denote by $\pi_X$ and $\pi_Y$ the corresponding
projections of $X\times Y$ to $X$ and $Y$ respectively.

By Mukai (\cite{Muk2}, Theorem 1.5) the
corresponding quasi-universal sheaf $\E$ on $X\times Y$ and defined by
this sheaf algebraic cycle
\begin{equation}
Z_{\E}=\pi^\ast_X(\sqrt{td_X})\cdot ch(\E)\cdot
(\pi^\ast_Y\sqrt{td_Y})/\sigma (\E)
\label{Mukaicycle}
\end{equation}
(see (\cite{Muk2}, Theorem 1.5) for details) define the
isomorphism of the full cohomology groups and the corresponding
Hodge structures
\begin{equation}
\label{Mukaicycleisom}
f_{Z_\E}:H^\ast(X,\bq)\to H^\ast (Y,\bq),\ \
t\mapsto \pi_{Y,\ast}(Z_\E\cdot\pi^\ast_X t).
\end{equation}
Moreover, according to Mukai,
it defines the isomorphism of lattices (or isometry)
$$
f_{Z_\E}:v^\perp \to H^4(Y,\bz)\oplus H^2(Y,\bz)
$$
where $f_{Z_\E}(v)=w\in H^4(Y,\bz)$ is the fundamental cocycle, and
the orthogonal complement $v^\perp$ is taken in the Mukai lattice
$\widetilde{H}(X,\bz)$. It follows \eqref{Mukai1} which we used in
Sect. \ref{secpicard1}.

In particular, taking the composition of $f_{Z_\E}$ with the
projection $\pi:H^4(Y,\bz)$ $\oplus H^2(Y,\bz)\to H^2(Y,\bz)$,
we obtain an embedding of lattices
$$
f_{Z_\E}:H^\perp_{H^2(X,\bz)}\to H^2(Y,\bz)
$$
which can be extended to the isometry
\begin{equation}
\label{isometryH2}
\widetilde{f}_{Z_\E}:H^2(X,\bq)\to H^2(Y,\bq)
\end{equation}
of the quadratic forms over $\bq$ by Witt's Theorem.

If $H^2=0$, this extension is unique.

If $H^2\not=0$, there are two such extensions different by $\pm 1$ on
$\bz H$. Let us agree to take
\begin{equation}
\label{isometryH2h}
\widetilde{f}_{Z_\E}(\wH)=c\,h
\end{equation}
where $h$ is defined in \eqref{definitionofh} and we use Proposition
\ref{propositionperY} which relates periods of $X$ and $Y$.

The defined Hodge isometry \eqref{isometryH2}
can be considered as a little change of the Mukai's algebraic
cycle \eqref{Mukaicycle} to get an isometry in $H^2$. Clearly, it is
also defined by some algebraic cycle because changes the Mukai
isomorphism \eqref{Mukaicycleisom} only in the algebraic part.

\medskip

By Proposition \ref{propositionperY}, we obtain that the isomorphism
$\widetilde{f}_\E$ is given by the embeddings
\begin{equation}
\label{corrXY} \wH^\perp \subset h^\perp=[\wH^\perp,2abct^\ast
(\wH) ],\ \bz \wH\subset \bz h,\ \wH=ch, H^{2,0}(X)=H^{2,0}(Y).
\end{equation}
This identifies quadratic forms $H^2(X,\bq)=H^2(Y,\bq)$ over $\bq$,
and the lattices $H^2(X,\bz)$, $H^2(Y,\bz)$ as its two sublattices.

Let
$$
O(H^2(X,\bq))_0=\{f\in O(H^2(X,\bq))\ |\ f|T(X)=\pm 1\}\cong
$$
$$
\cong O(N(X)\otimes \bq)\times \{\pm 1_{T(X)}\},
$$
and
$$
O(H^2(X,\bz))_0=O(H^2(X,\bz))\cap O(H^2(X,\bq))_0.
$$
By Global Torelli Theorem for K3 surfaces \cite{PS}, we obtain, at once:

\begin{proposition}
\label{propGlToraction}
If a K3 surface $X$ is general for its Picard lattice,
then $Y=M_X(r,H,s)\cong X$ if and only if
there exists an automorphism
$\phi(r,H,s)\in O(H^2(X,\bq)_0)$ such that $\phi(H^2(X,\bz))=H^2(Y,\bz)$.
\end{proposition}

If $Y\cong X$, then we can give the definition.

\begin{definition}
\label{definitionaction} If $Y=M_X(r,H,s)\cong X$ and $X$ is general for its
Picard lattice, then the isomorphism
of Proposition \ref{propGlToraction}
$$
\phi (r,H,s) \mod O(H^2(X,\bz))_0
\in O(H^2(X,\bq))_0/O(H^2(X,\bz))_0
$$
is called the {\bf action on $H^2(X,\bq)$} of the correspondence
of a general (for its Picard lattice)
K3 surface $X$ with itself via moduli of sheaves $Y=M_X(r,H,s)$ on $X$ with
the primitive isotropic Mukai vector $v=(r,H,s)$.
\end{definition}

By Global Torelli Theorem for K3 surfaces \cite{PS},
the group $O(H^2(X,\bz))_0\mod {\pm 1}$
can be considered as generated by correspondences defined by graphs of
automorphisms of $X$ and by the reflections
$s_\delta:x\mapsto x+(x\cdot \delta)\delta$, $x\in H^2(X,\bz)$,
in elements $\delta\in N(X)$ with $\delta^2=-2$. By Riemann--Roch Theorem
for K3 surfaces, $\pm \delta$ contains an effective curve $E$.
If $\Delta\subset X\times X$ is the diagonal, the effective
2-dimensional algebraic cycle $\Delta+E\times E\subset X\times X$
acts as the reflection $s_\delta$ in $H^2(X,\bz)$
(the author knows this from Mukai \cite{Muk5}). Thus, considering
actions of correspondences modulo $O(H^2(X,\bz))\mod {\pm 1}$ is very
natural.

\medskip

Let us consider Tyurin's isomorphism \ref{Tyurinisom} defined by the
Mukai vector $v=(\pm H^2/2,H,\pm 1)$ where $H\in N(X)$ has $H^2\not=0$
and $\pm H^2>0$. Assume that $\wH=H/d$ is primitive.
Then
$M_X(\pm H^2/2,H,\pm 1)\cong M_X(\pm \wH^2/2 \ \wH,\pm 1)$,
and we can assume that $\wH$ is primitive.

Then $c=1$, $a=\pm \wH^2/2$ and $b=\pm 1$, $m(a,b)\equiv -1\mod 2a$,
$h=\wH$.

We have
$$
H^2(X,\bz)=[\bz \wH, \wH^\perp, \wH+t^\ast (\wH) ],
$$
$$
H^2(Y,\bz)=[\bz \wH, \wH^\perp, \wH-t^\ast (\wH) ],
$$
Then the reflection $s_{\wH}$ with respect to $\wH$,
$$
s_{\wH}(x)=x-\frac{2(x \cdot \wH)\wH}{\wH^2},\ \ x \in H^2(X,\bq),
$$
belongs to $O(H^2(X,\bq))_0$ and $s_{\wH}(H^2(X,\bz))=H^2(Y,\bz)$.

Thus, we obtain

\begin{proposition} For a K3 surface $X$ and
$H\in N(X)$ with $\pm H^2>0$, the Tyurin isomorphism
$$
M_X(\pm H^2,H,\pm 1)\cong X
$$
defines the correspondence of $X$ with itself with the action
$$
s_H\mod O(H^2(X,\bz)_0)
$$
where $s_H$ is the reflection in the element $H$.

By classical and well-known results, their compositions generate
the full group $O(H^2(X,\bq))_0\mod {\pm 1}$.
\end{proposition}

\subsection{The general problem of classification of correspondences of a
K3 surface with itself via moduli of sheaves}
\label{subsecprobclass}

We will need some notations. For a sublattice $N\subset N(X)$, we introduce
$$
O(N\otimes \bq)_0=\{f\in O(H^2(X,\bq))\ |\ f|N^\perp_{H^2(X,\bz)}=\pm 1\}
$$
and
$$
O(N)_0=O(H^2(X,\bz))\cap O(N\otimes \bq)_0.
$$

Let $X$ be a general (for its Picard lattice) K3 surface $X$ and $N(X)$
its Picard lattice. The problem of classification of correspondences of $X$
with itself via moduli of sheaves consists of the following problems
(from the author's point of view now):

\medskip

{\it (1) Find all primitive isotropic Mukai vectors $(r,H,s)$ on $X$
such that $Y=M_X(r,H,s)\cong X$.}

\medskip

{\it (2) For a primitive isotropic Mukai vector
$(r,H,s)$ from (1), find all critical polarized Picard sublattices
$H\in N(r,H,s)\subset N(X)$.}

\noindent
For each of them, the corresponding
$\phi(r,H,s)$ from Definition \ref{definitionaction}
can be taken from $O(N(r,H,s)\otimes \bq)_0$ (we denote it as
$\phi_{N(r,H,s)}$, and it looks like
a reflection with respect to $N(r,H,s)$. For two critical polarized
Picard sublattices
$H\in N(r,H,s)$ and $H\in N^\prime(r,H,s)$,
the automorphisms $\phi_{N(r,H,s)}$ and $\phi_{N^\prime(r,H,s)}$
are different by an automorphism from $O(H^2(X,\bz))_0$.

\medskip

{\it (3) The structures (1) and (2) are important}
because for two primitive isotropic
Mukai vectors $(r,H,s)$ and $(r^\prime,H^\prime,s^\prime)$ from (1) and
two their critical polarized Picard sublattices $H\in N(r,H,s)$ and
$H^\prime \in N(r^\prime, H^\prime, s^\prime)$, the isomorphism between
$M_X(r,H,s)$ and $M_X(r^\prime, H^\prime, s^\prime)$ which is defined by
$$
\phi(r^\prime,H^\prime, s^\prime)\phi (r,H,s)^{-1},
$$
comes from K3 surfaces with the Picard sublattice
$N(r,H,s)+N(r^\prime, H^\prime, s^\prime)\subset N(X)$, and it can be
considered as a natural isomorphism between these moduli.

\medskip

(4) All these generators $\phi_{N(r,H,s)}\mod O(N(r,H,s))_0$
can be considered as natural generators for correspondences of $X$ with
itself via moduli of sheaves, together with automorphisms of $X$ and
reflections $s_\delta$, $\delta\in N(X)$ and $\delta^2=-2$.
They and their relations are the natural subject to study.

\medskip

For $\rho (X)=1,\ 2$, problems (1)---(4) are solved.
See  Sections \ref{secpicard1} and \ref{secpicard2}. Results of
Sect. \ref{secpicard3} show that these problems are very non-trivial for
$\rho (X)\ge 3$.

\medskip

As an example, let us take a general K3 surface $X$ with the Picard lattice
$N(X)=S$ of Theorem \ref{theoremrho=3} of the rank three (or any other
Picard lattice of the rank three which satisfies Theorem
\ref{theoremuniqtrper}). Let $v=(r,H,s)$ be a primitive isotropic Mukai vector
on $X$. Then $Y=M_X(r,H,s)\cong X$ if and only if $(\gamma,c)=1$. Then we
have three cases:

(a) If $c=1$ and either $a_1=1$ or $b_1=\pm 1$ (Tyurin's case), then the
critical sublattice is $N(v)=\bz \wH$, it  has the rank one and is unique.
The corresponding $\phi_{N(v)}=s_H\mod O(H^2(X,\bz))_0$.

(b) If $v=(r,H,s)$ is different from (a), but the critical sublattice
$N(v)$ has the rank two (the divisorial case), then all critical
sublattices $N(v)$ are generated by $H$ and $h_1\in [H, a_1cN(X)]$ with
$h_1^2=\pm 2a_1c$ or $h_1\in [H,b_1cN(X)]$ with $h_1^2=\pm 2b_1c$.
(See theorems of Sect. \ref{secpicard2}). All these $N(v)$
give automorphisms $\phi_{N(v)}$ which are different by elements from
$O(H^2(X,\bz))_0$.

(c) If $v=(r,H,s)$ is different from (a) and (b), then the critical
sublattice $N(v)=N(X)$ has the rank three. These cases really happen
by Theorem \ref{theoremrho=3}. We obtain
$\phi_{N(v)}\mod O(H^2(X,\bz))_0$.

\medskip

Any two $v_1$, $v_2$ satisfying one of these conditions (a), (b) or (c),
and corresponding their critical sublattices $N(v_1)$, $N(v_2)$ generate
natural isomorphisms $\phi_{N(v_2)}\phi_{N(v_1)}^{-1}$ between corresponding
moduli of sheaves over $X$ (all of them are isomorphic to $X$)
which are specializations of the corresponding isomorphisms from
the Picard sublattice $N(v_1)+N(v_2)\subset N(X)$.

\medskip

A reader can see that our general idea is that a very complicated
structure of correspondences of a general (for its Picard lattice)
K3 surface $X$ with itself via moduli of sheaves is hidden inside
of the abstract Picard lattice $N(X)$, and we try to recover this
structure. This should lead to some non-trivial constructions
related to the abstract Picard lattice $N(X)$ and more closely
relate it to geometry of the K3 surface.


\

\

V.V. Nikulin \par Deptm. of Pure Mathem. The University of
Liverpool, Liverpool\par L69 3BX, UK; \vskip1pt Steklov
Mathematical Institute,\par ul. Gubkina 8, Moscow 117966, GSP-1,
Russia

vnikulin@liv.ac.uk \ \ vvnikulin@list.ru

\end{document}